# DESIGN OF COOPERATIVE PROCESSES IN A CUSTOMER-SUPPLIER RELATIONSHIP: AN APPROACH BASED ON SIMULATION AND DECISION THEORY

**François GALASSO**

LAAS-CNRS, Université de TOULOUSE
7 avenue du Colonel ROCHE
31077 Toulouse Cedex 4
galasso@univ-tlse2.fr

**Caroline THIERRY**

Université de TOULOUSE, IRIT
5 allées Antonio MACHADO
31058 TOULOUSE
thierry@univ-tlse2.fr

**ABSTRACT:** *Performance improvement in supply chains, taking into account customer demand in the tactical planning process is essential. It is more and more difficult for the customers to ensure a certain level of demand over a medium term horizon as their own customers ask them for personalisation and fast adaptation. It is thus necessary to develop methods and decision support systems to reconcile the order and book processes. In this context, this paper intends firstly to relate decision under uncertainty and the industrial point of view based on the notion of risk management. This serves as a basis for the definition of an approach based on simulation and decision theory that is dedicated to the design of cooperative processes in a customer-supplier relationship. This approach includes the evaluation, in terms of risk, of different cooperative processes using a simulation-dedicated tool. The evaluation process is based on an exploitation of decision theory concepts and methods. The implementation of the approach is illustrated on an academic example typical of the aeronautics supply chain.*

**KEYWORDS:** *supply chain, simulation, cooperation, decision theory, risk*

## 1 INTRODUCTION

Supply chain management emphasises the necessity to establish cooperative processes that rationalize or integrate the forecasting and management of demand, reconcile the order and book processes, and mitigate risks.

These cooperative processes are often characterised by a set of point-to-point (customer/supplier) relationships with partial information sharing (Galasso *et al.*, 2006). In this context, at each level of the supply chain, a good understanding of the customer demand is a key parameter for the efficiency of the internal processes and the upstream supply chain (Bartezzaghi and Verganti, 1995). However, due to a substantial difference among the supply chain actors in terms of maturity regarding their use of enterprise systems, it is more or less difficult to implement cooperative processes for the different participating companies. Indeed, while large companies have the capability of using and managing efficient cooperative tools, Small and Medium Enterprises (SMEs) suffer from a partial vision of the supply chain and have difficulties to analyse the uncertain information communicated from customers.

This paper aims at providing suppliers (e.g. in aeronautics) with a cooperation support that takes advantage of the information provided by customers in a cooperative perspective even if this information is uncertain. Thus, we propose a decision and cooperation support approach based on a simulation of planning processes in the point-to-point supply chain relationship.

More precisely, we are concerned with the joint evaluation of the impact of the customer's supply management process and the supplier's demand management and planning processes.

After discussing the state of the art (cf. Section 2) on cooperation in supply chain management and Supply Chain Risk Management, we introduce the context and the related challenges (cf. Section 3). Then, section 4, describes the approach based on simulation and decision theory proposed to evaluate the risks pertained to the choice strategies for demand management (supplier) and supply management (customer). At last, the proposed methodology is implemented on an illustrative example (cf. Section 5).

## 2 STATE OF THE ART

In this section, the state of the art regarding two main points of interest is given. The first issue refers to decision making under uncertainty in connexion with the Supply Chain Risk Management. The second one investigates the problematics of cooperation within the supply chain.

### 2.1 Decision under uncertainty and Supply Chain Risk Management (SCRM)

In the industrial context, the concept of decision under uncertainty is generally not explicitly addressed but the concept of risk management is prominent. It is undeniable that the concepts of "*uncertainty*" and "*risk*" are linked even if it is sometimes difficult to perceive this link. Risk management, particularly in the field of supply chain management, turns out to be an important industrial challenge. Supply Chain Risk Management (SCRM) is the "management of external risks and supply chain risks through a coordinated approach between the supply chain partners in order to reduce supply chain vulnerability as a whole" (Christopher, 2003). So far, there is still a "lack of industrial experience and academic research for supply chain risk management" as identified by (Ziegenbein and Nienhaus, 2004) even if, since 2004, it has been an increasing number of publications in this field. More specifically, the question of risk management related to the use of Advanced Planning Systems has to be studied (Ritchie *et al.*, 2004). The academic community paid a lot of attention to the clarification of definition, taxonomies and models linked to the SCRM ((Brindley, 2004), (Tang, 2006)). Holton (2004) defines risk as the combination of two mains elements: the exposure and the uncertainty. Thus, he defines risk as the "exposition to a proposition (i.e. a fact) that one is uncertain".

However, from the viewpoint of decision theory, the distinction between "decision under risk" and "decision under uncertainty" is well established according to the knowledge of the state of the nature: the term decision under risk is used if objective probabilities are associated to the occurrences and if not, the term decision under uncertainty is used (Lang, 2003). The latter corresponds to the situation of imperfect knowledge. The imperfection of the knowledge of a system can be due to the flexibility inherent to the knowledge or due to the acquisition of such knowledge. Among these imperfections, Bouchon-Meunier (1995) synthesizes the distinction (Dubois and Prades, 1988): uncertainty; imprecision and incompleteness. Uncertainty refers to the "doubt about the validity of the knowledge", which refers to the fact of being unsure whether a proposition is true or not (for example: "I believe but I am not sure"…). Imprecision concerns "the difficulty to express knowledge". Indeed, it can be knowledge expressed in natural language in vague way (for example: "it is important"…) or quantitative knowledge not precisely known because of, for example, imprecise measurement ("this value lies between x and y" or "this value can be x, y or z"). Incompleteness refers to "the lack of knowledge or partial knowledge about some characteristics of the considered system".

A lot of criteria can be used in order to finely classify the different kinds of uncertainties (Teixidor, 2006). Bräutigam *et al.*, (2003) distinguish between two main kinds of uncertainties: endogenous uncertainty (specific to the studied company or system) and exogenous uncertainties (external to the studied company or system). More precisely in the field of Supply Chain Management, Ritchie *et al.*, (2004) propose a contingency framework over 4 dimensions: the environment characteristics, the supply chain context, the decisional system (decision level, type of decision, information availability,…), the human and its behaviour in presence of risk.

Regarding the production planning models under uncertainty, Mula *et al.*, (2006) have recently proposed a complete state of the art. In this review, the authors distinguish conceptual models; analytical models; artificial intelligence models and simulation models in order to deal with risk management issues. In the last category, the model proposed by Rota *et al.*, (2002) can be pointed out as it is close to our approach embedding an analytical model in a simulation framework. Indeed, it is one of the first attempts made in order to evaluate the interest of taking into account forecasts in the planning process while software such as Advanced Planning Systems (APS) just began to be implemented. Nowadays, considering the spreading out of the use of such tools, practitioners aim at quantifying the risk inherent to the planning process with an APS in the supply chain context (Ritchie *et al.*, 2004). In that sense, Génin *et al.*, (2007) propose, for example, an approach that provides a robust planning with an APS. Beyond the planning process in itself, it becomes more and more important to assist industrial practitioners in defining demand management in order to deal with uncertainty while maximising the potential use of the planning tools.

## 2.2 Cooperation in Supply chains

The implementation of cooperative processes for supply chain management is a central concern for practitioners and researchers. This awareness is linked, in particular, to the Bullwhip effect whose influence has been clearly shown and studied (Lee *et al.*, 1997; Moyaux, 2004).
Recently, many organizations have emerged to encourage trading partners to establish cooperative interactions (that rationalize or integrate their demand forecasting/management, and reconcile the order-book processes) and to provide standards (that could support cooperative processes): RosettaNet (Rosetta, 2007), Voluntary Inter-industry Commerce Standards Association (Vics, 2007), ODETTE (Odette, 2007), etc. On the other hand, McCarthy and Golicic (2002) consider that the cooperative process brought by the CPFR (Collaborative Planning, Forecasting and Replenishment) model is too detailed. They suggest instead that companies should plan regular meetings to discuss the forecast with the other supply chain partners so as to develop shared forecast.

In the same way, many recent research papers are devoted to cooperation in the context of supply chain management. Under the heading of cooperation, authors list several aspects. One of these aspects on which we focus in this paper, is cooperation through information sharing. Using Huang *et al.* (2003) literature review, we can distinguish between different classes of information that play a role in the information sharing literature: (i) product information, (ii) process information, (iii) lead time, (iv) cost, (v) quality information, (vi) resource information, (vii) order and inventory information, (viii) Planning (forecast) information (Lapide, 2001; Moyaux, 2004). Another aspect of cooperation concerns that extend information sharing to collaborative forecasting and planning systems (Dudek and Stadtler, 2005; Shirodkar and Kempf, 2006). In this paper, we will focus on information sharing and more precisely sharing information concerning planning (forecast).

Nevertheless little attention has been paid to the risk evaluation of new collaborative processes (Småros, 2005, Brindley, 2004, Tang et al 2006). This is also true when planning processes under uncertainty are concerned (Mula *et al.,* 2006) even if Rota *et al.,* (2002) introduced the problem of managing tactical planning with an APS and Génin *et al.,* (2007) studied its robustness.
Thus, this paper is focused on risk evaluation of cooperative planning processes within a customer-supplier relationship and thus, a decision and cooperation support tool for dealing with uncertainty is proposed.

## 3 CONTEXT AND CHALLENGES

It has been stressed in section 2 how building cooperative processes is of major importance. The main concern regarding cooperation in a supplier-customer relationship context is not to argue about the interest of the cooperation but to define and support the cooperative process. Along this line, two viewpoints can be adopted in order to manage the supply chain: supply chains can be managed, on the one hand, as a single entity through a dominant member […] and on the other hand, through a system of partnerships requiring well-developed cooperation and coordination (Ganeshan *et al.,* 1999). In this paper, the second viewpoint is adopted. Moreover, in this system, we focus on a partnership between a customer and a supplier (SME) (cf. Figure 1.). In the considered relationship, the two actors are juridically independent and in charge of their own planning processes (embedding their own suppliers and subcontractors). Nevertheless, a global distributed planning process exists which includes the individual planning and information exchange. This planning process is the result of a cooperative design involving the two partners.

More precisely, on the customer side, the demand management process is studied. The customer communicates a demand plan to the supplier with a given periodicity.

This plan is established thanks to the customer planning process in which a frozen horizon is considered (within this frozen horizon no decision can be revised). *Firm demands* related to the period close to the present time are communicated to the supplier within this frozen horizon. They are defined on a given time horizon, called firm horizon (*FH*).

Beyond this horizon, decisions can be revised within a given interval. This interval is part of the cooperation partnership between the supplier and the customer. We call "forecast" or "flexible" demands the pair (forecast value, flexibility level) which is communicated to the supplier. The flexibility is expressed in term of percentage of variation around the forecast value. The minimum and maximum values of the flexibility interval will be called "flexibility bounds" here after. This flexible demand is defined on a given time horizon, called flexible horizon (*LH*) which is considered as part of the cooperation framework between the customer and the supplier.

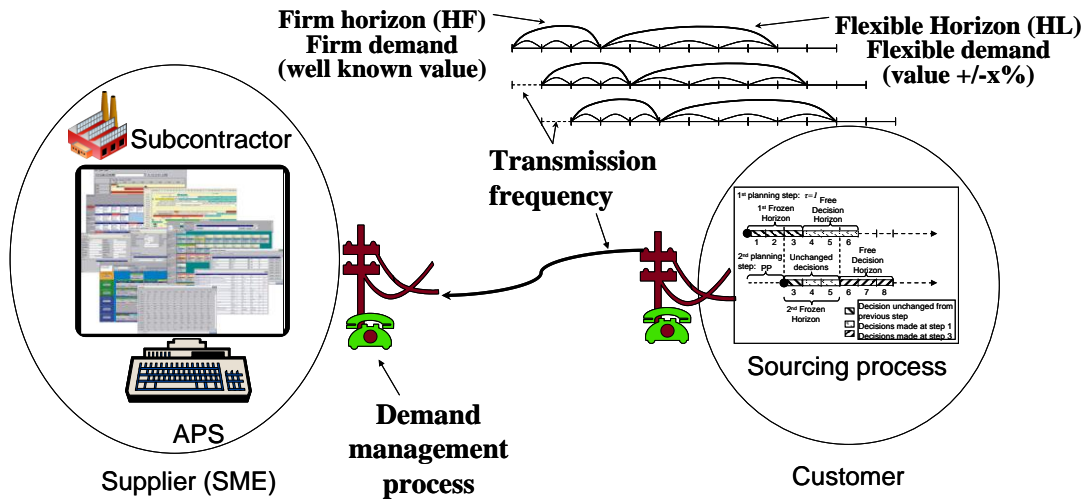

Figure 1. Study positioning

Moreover, on the supplier side, the planning process at a given time point is performed through the optimisation procedure of an Advanced Planning System, which is not the main object of this study. The APS computes deterministic data thus the supplier has to pre-compute the flexible demands communicated by the customer. Different types of behaviours are studied according to the degree of knowledge of the supplier on his customer's behaviour (for example, propensity to overestimate or to underestimate).

Adopting the supplier point of view, uncertainty is mainly due to the trends of the forecasted requirements of the customer embedded in the demand plan communicated to the supplier. The uncertainty about the trends can be due to:
- the difficulty to forecast the global market evolution that could be increasing, decreasing, or having a temporary variation (i.e. a peak of demand),
- the distribution of the customer's production on the considered horizon

These requirements can be considered as *imprecise*. In the illustrative example developed in section 5, it is assumed that the customer knows that there will be a peak in the demand. However, some *uncertainty* remains about the height of the peak. Moreover, the customer's requirements at a given date are given between two bounds: the lower and the upper bounds of the flexibility interval.

Furthermore, the demands are communicated by the customer over a given time horizon. Beyond this horizon, no demand is expressed and thus, the knowledge about the forecasted demand is *incomplete*.

Adopting the customer point of view, there is no information regarding the demand management strategy of the supplier. Thus, it remains uncertain information for the customer.

In this context, the challenge is to support the decision makers in order to set up the supplier's demand management behaviour and the customer's forecast transmission behaviour within a cooperative supply chain planning process. Thus, an approach based on simulation and decision theory enabling a risk evaluation of the demand management processes according to different scenarios is proposed.

# 4 DECISION AND COOPERATION SUPPORT UNDER UNCERTAINTY

In this section, a decision and cooperation support approach for designing the strategies of demand management (on the supplier side) and requirement management (on the customer side) is presented. This approach is presented in section 4.1. It uses a simulation tool detailed in section 4.2 which embeds a model for the behaviour of both actors of the considered relationship.

## 4.1 Decision and cooperation support approach using simulation

As one of our main goals is to create reliable partnerships, it is necessary to establish a discussion between the decision makers of the different entities. Thus, an implementation methodology has been adapted which can be set up in five steps (Lamothe *et al.,* 2008).

Step 1: Problem and system definition

Step 2: Design of experiment

Step 3: Simulations (performed over a given time horizon defined by both actors) and computation of the indicators

Step 4: Risk evaluation: the risks of the different cooperation strategies are expressed for both actors of the chain (The word "risk" being used here in its industrial acceptance)

Step 5: Design of cooperative processes

The different steps of the proposed approach are detailed hereinafter.

### 4.1.1 Problem and system definition (step 1)

This step is a viewpoint confrontation process and is led by an external organizer (for example the tool designer) (Thierry *et al.,* 2006) with:
- The presentation of the actors concerned and the expression of the context in which the decision takes place.
- The definition of the fundamentals of the potential cooperative planning strategies to be evaluated,
- The definition of the risks to be evaluated (for example the risk of a supplier strategy)
- The choice of the indicators enabling risk evaluation (for example the global gain of the supplier)

In the context considered in this study, we consider a relationship including a customer and a supplier. Both actors have to determine they behaviours (internal strategies) to design a common cooperative strategy.
The main problem of the supplier is to choose a planning strategy concerning the demand management in order to take into account the demand communicated by the customer in its planning process.
Regarding the customer's side, the supply management process is considered. Within this process, an important decisional lever is the length of the firm and flexible horizon. Through this lever, the supplier has more or less visibility on the demand and thus more or less time to react and adapt its production process.
At the supplier level, the definition of a cost model (a cost being associated to each parameter of the model given in section 4.2.2) enables the calculation of *the global gain* obtained by the use of each strategy regarding each scenario. This gain can be considered as representative, at an aggregated level, of the combination of all indicators to evaluate the risks associated to the planning policies that he envisaged. Nevertheless, as the problem cannot be totally defined at once, other indicators have to be computed (stock levels, service level, production cost, …). At the customer level, the *cost of backorders,* for example, can be considered as pertinent.

### 4.1.2 Design of experiment (steps 2)

In this step, for each actor of the supply chain, the set of values of the strategies parameters are defined. Moreover, different scenarios are defined by a combination of values of parameters which are the uncontrolled variables of the considered actor. An experiment is defined by the combination of all these parameters.

### 4.1.3 Simulations and computation of the indicators (steps 3)

In this step a dedicated simulation tool is used to run the design of experiment and the indicators defined in step 1 are computed.

### 4.1.4 Risk evaluation (steps 4)

This risk evaluation step is at the heart of the approach. It is based on an implementation of the criteria commonly used decision theory. In order to engage a cooperative process, it is necessary to consider the objectives of both actors of the supply chain simultaneously. To perform this multi-actor decision making process, we propose a cooperative dashboard (an example is given on Figure 2).

#### 4.1.4.1 Cooperative decision support dashboard

The cooperative decision support is divided into two sides: the supplier side and the customer's side. These two sides are shared by the two actors.
For each side we propose a so-called risk diagram (the design of this risk diagram is detailed in section 4.1.3.2.) which is the central decision support for the planning strategies choice. Moreover a regret table is proposed (cf. section 4.1.3.3.) to the decision makers which enables the proposed strategies to be situated within the set of potential strategies. Then, we associate to each selected strategy a set of other indicators (cf. section 4.1.3.4) measuring inventory, production, purchasing levels…

#### 4.1.4.2 Risk diagram

From an actor point of view, the best policy can be different depending on the considered scenario. Thus, it is necessary to compare each strategy considering the whole set of scenarios. In such a context, such a comparison is possible using a decision criterion in order to aggregate the indicators obtained for each scenario. In the framework of the problem under study, it is hardly possible to associate probabilities to the occurrence of each scenario. Thus, the evaluation can be done through the use of several decision criteria (which may lead to different results) based on the gain or the costs obtained after the simulation of each scenario: Laplace's criterion (average evaluation), Wald's criterion (pessimistic evaluation), Hurwicz's criterion (weighted sum of pessimistic and optimistic evaluation), Savage's criterion (minimising the maximum regret), etc. The results given by the different criteria can be gathered into a risk diagram on which managers in charge of the planning process can base their decision making (Mahmoudi, 2006). A general diagram is presented and detailed in Figure 3.

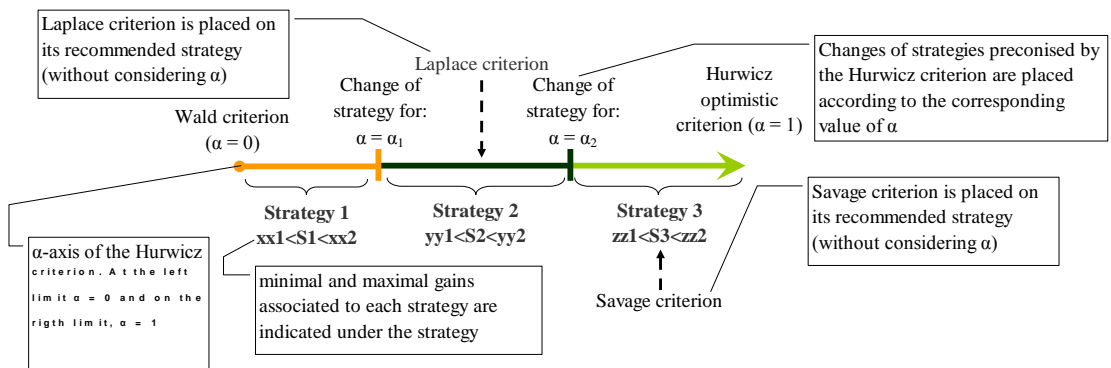

Figure 2. General risk diagram

In this diagram, the demand management strategies are positioned regarding the risk-prone attitude of the decision maker: these strategies are thus positioned on an axis corresponding to the values of α between 0 and 1 and denoted α-axis. The evolution of the value of this criterion as a function of α for each strategy is represented on a curve following the formula of the Hurwicz criterion: $H_S(α) = (1-α)\, m_S + α\, M_S$ (with $m_S$ the minimal gain and $M_S$ the maximal gain obtained applying the strategy $S$). From this curve, the values of $α_i$ indicating a change in the proposed strategy can be determined. Then, the strategies are specified on the diagram. For each strategy, the associated minimal and maximal gains or costs (according to the selected indicator as, for example, the global gain or the backorder costs) are given. Furthermore, if the represented strategies are evaluated by mean of other criteria (Laplace or Savage), these criteria are attached to the relevant strategy (without considering the value of α).

#### 4.1.4.3 Regret table

The risk diagram is not sufficient to have en exhaustive comparison of strategies. Thus, the purpose of the regret table is to give an indication about the risk taken when using a strategy instead of another one. The regret of using each strategy regarding the others is given to each actor of the supply chain. This regret is calculated as the difference between the gain obtained with the strategy which could be used and a reference strategy. For each pair of strategies, the minimal and the maximal regrets are given as depicted in Table 1.

|                 | Regret using K1 | Regret using K2 | Regret using K3 |
|-----------------|-----------------|-----------------|-----------------|
| Min regret / K1 | 0               | -240            | -6 400          |
| Max regret / K1 | 0               | 300             | 900             |
| Min regret / K2 | -300            | 0               | -440            |
| Max regret / K2 | 240             | 0               | 600             |
| Min regret / K3 | -900            | -600            | 0               |
| Max regret / K3 | 6 400           | 440             | 0               |

Table 1. Illustrative regrets table

Table 1 considers generic strategies K1, K2 and K3 for the sake of illustration. Table 1 provides a skew-symmetric matrix as, for example, the regret of using K1 instead of K2 is the opposite of the regret of using K2 instead of K1.

#### 4.1.4.4 Other dashboard indicators

Beyond the risk diagram and the regret table, the dashboard provides the decision makers with information regarding the production conditions. On the supplier's side, the global gains, inventory costs and production costs are given. This information is completed with specific inventory and production evolution in the cases giving the highest production and inventory costs. On the customer's side, the evolution of backorders is particularly studied.

| Supplier Dashboard | Customer Dashboard |
|---|---|
| ### Risk Diagram | ### Risk Diagram |

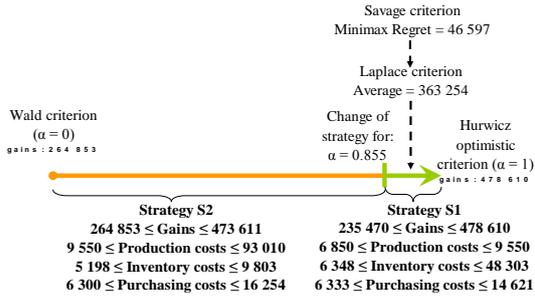

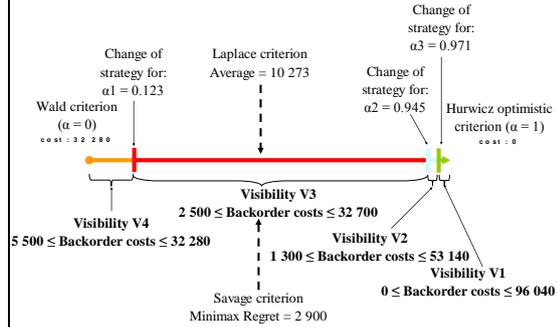

### Regret Table

|  | Regret using S1 | Regret using S2 |
|---|---|---|
| Min regret / S1 | 0 | -46 597 |
| Max regret / S1 | 0 | 73 034 |
| Min regret / S2 | -73 034 | 0 |
| Max regret / S2 | 46 597 | 0 |

### Regret table

|  | Regret using V1 | Regret using V2 | Regret using V3 | Regret using V4 |
|---|---|---|---|---|
| Min regret / V1 | 0 | -43 240 | -63 400 | -63 760 |
| Max regret / V1 | 0 | 1 300 | 2 900 | 5 900 |
| Min regret / V2 | -1 300 | 0 | -20 440 | -20 860 |
| Max regret / V2 | 43 240 | 0 | 1 600 | 4 600 |
| Min regret / V3 | -2 900 | -1 600 | 0 | -420 |
| Max regret / V3 | 63 400 | 20 440 | 0 | 3 000 |
| Min regret / V4 | -5 900 | -4 600 | -3 000 | 0 |
| Max regret / V4 | 63 760 | 20 860 | 420 | 0 |

### Comparison of gains

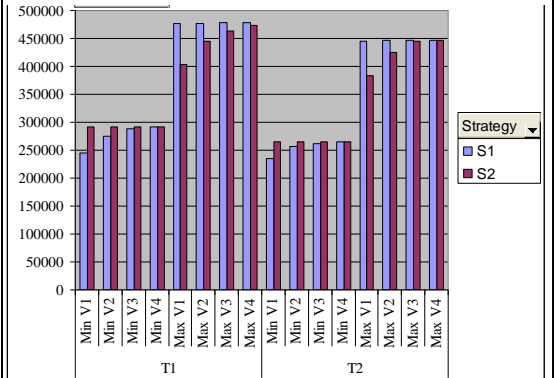

### Comparison of backorder costs

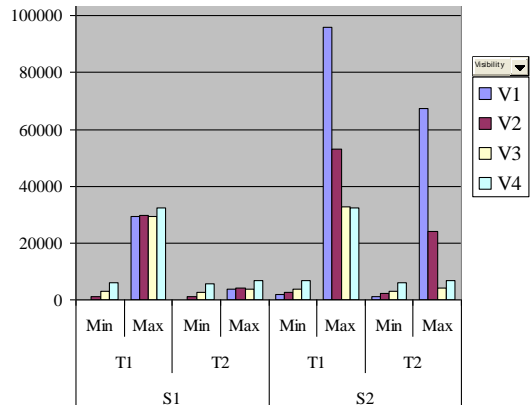

### Comparison of inventory costs

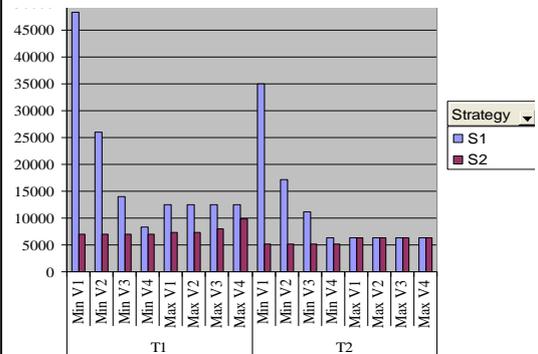

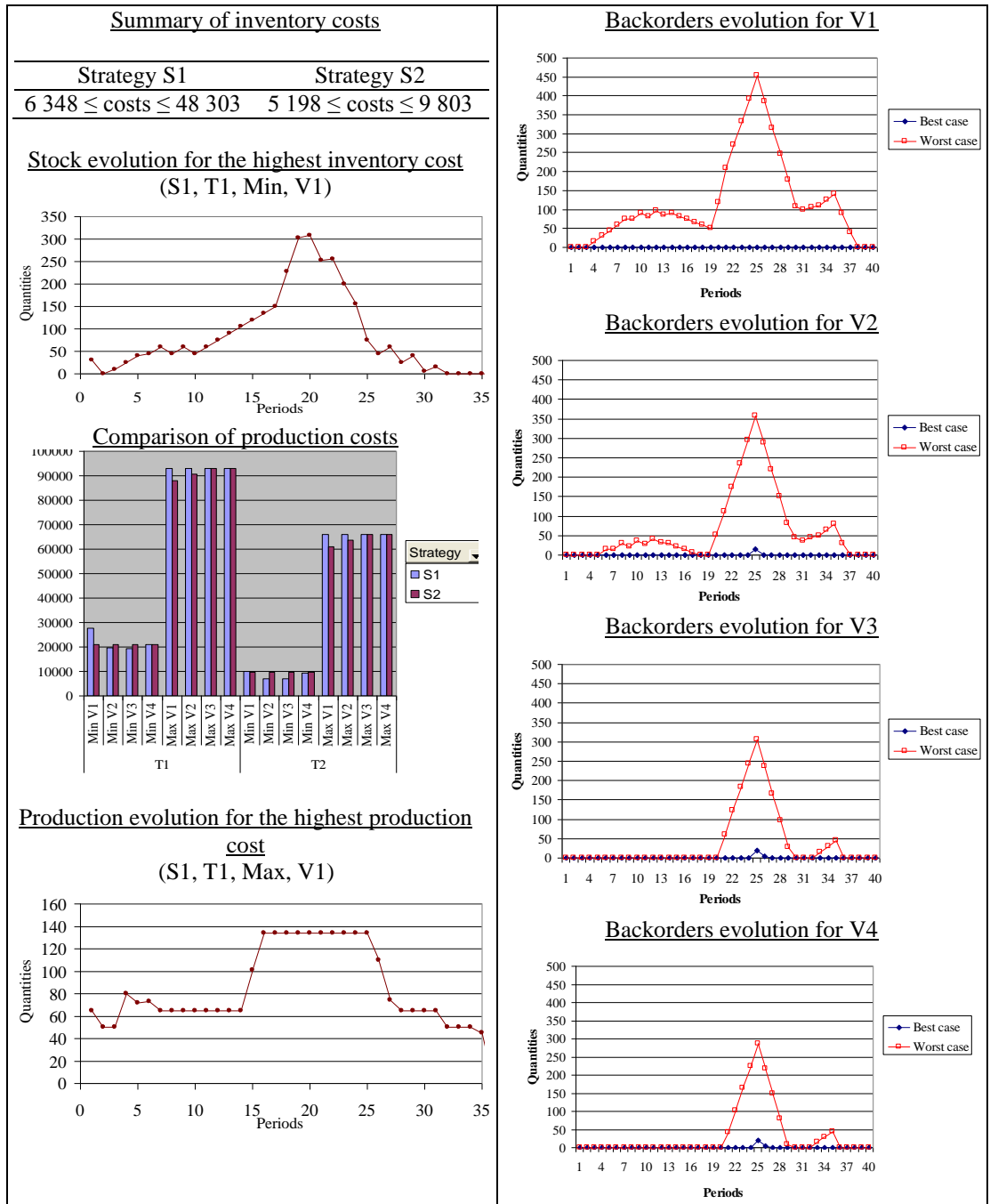

Figure 3. Example of cooperative dashboard

### 4.1.5 Design of cooperative processes (step 5)

Using the dashboard of the preceding step, the decision makers are able to analyze the simulation using the proposed dashboard with:
- a risk analysis from the supplier's point of view,
- a risk analysis from the customer's point of view,
- a risk analysis from a cooperative point of view

This analyse may lead either
- to a common identification of a need for a second run of the decision and cooperation support approach with the redefinition of the criteria of the risk evaluation, of the risks themselves, of the design of experiment…

- or to the conclusion of the partnership contract pertaining to the choice of cooperative processes,

**4.2 The simulation tool**

In order to model the dynamic behaviour of both actors a dedicated simulation tool has been developed enabling the evaluation of:
- The supply management behaviour models of the customer including the computation of firm demand and forecasts and their transmission to the supplier,
- The behaviour models of the supplier embedding:
    o The management process of the demand
    o The planning process

The simulation of these behaviours relies on a fixed time step. This period corresponds to the replanning period.

**4.2.1 Model of the customer's behaviour (supply management process)**

A model enables a macroscopic point of view of the customer's behaviour concerning the supply management. This model simulates the evolution of the customer demand communicated to its supplier. Considering a given visibility level of the demand (size of the firm and the flexible horizons), this model computes the customer demand at each planning step.
The principle of this model is illustrated hereafter on an example (figure 4).

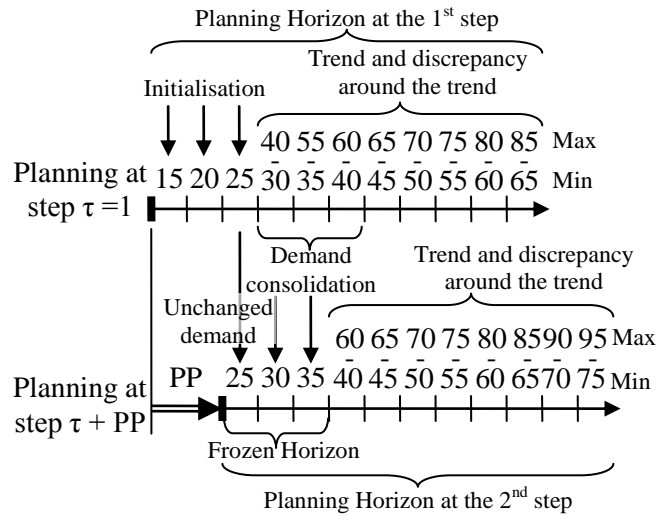

Figure 4. Customer's behaviour model

The procurement plan communicated to the supplier is established, on the flexible horizon, taking into account a trend and a discrepancy around this trend. The firm demand is calculated, on the firm horizon, according to the flexible demand established at the previous planning step and so-called hereinafter: the consolidation process of the demand. During the foremost planning step, the demand is initialised by the calculation of a flexible demand from the trend and the discrepancy over the whole planning horizon and then, the consolidation process is carried out over the firm horizon.

In the example depicted by Figure 4, the trend is linear and it grows-up to a 5 units production per period rate. The discrepancy is, in a simplified way, of +/- 5 units at each period. The modelled scenario is the one in which the customer overestimates the flexible demand. The firm demand is therefore calculated as equal to the lower bound of the communicated flexible demand at the previous planning step.

The customer demand is noted $D_{p,t}^{\tau}$. The discrepancy is modelled by an interval limited by the following bounds:
- $\underline{D}_{p,t}^{\tau}$, is the lower bound of the tolerated discrepancy over the flexible demand,

➤ $\overline{D}_{p,t}^{\tau}$, is the upper bound.

The demand expressed at each period are always within the interval defined by $\left[\underline{D}_{p,t}^{\tau}, \overline{D}_{p,t}^{\tau}\right]$ at each end-item $p$, period $t$ and planning step $\tau$. They are modelled as follows (1):

$$\begin{cases} D_{p,t}^{\tau}, & \forall p, \forall t \in FH^{\tau} \\ D_{p,t}^{\tau} \in \left[\underline{D}_{p,t}^{\tau}, \overline{D}_{p,t}^{\tau}\right] & \forall p, \forall t \in LH^{\tau} \end{cases} \quad (1)$$

The evolution of the demand between two successive steps is formalised by the following relations:

$$D_{p,t}^{\tau} = D_{p,t}^{\tau-PP} \quad \forall p \quad \forall t \in \{FH^{\tau-PP} \cap FH^{\tau}\} \quad (2)$$

$$D_{p,t}^{\tau} = g(\underline{D}_{p,t}^{\tau-PP}, \overline{D}_{p,t}^{\tau-PP}) \text{ with } D_{p,t}^{\tau} \in \left[\underline{D}_{p,t}^{\tau-PP}, \overline{D}_{p,t}^{\tau-PP}\right] \forall p \quad \forall t \in \{LH^{\tau-PP} \cap FH^{\tau}\} \quad (3)$$

$$\left[\underline{D}_{p,t}^{\tau}, \overline{D}_{p,t}^{\tau}\right] = \left[\underline{D}_{p,t}^{\tau-PP}, \overline{D}_{p,t}^{\tau-PP}\right] \forall p \quad \forall t \in \{LH^{\tau-PP} \cap LH^{\tau}\} \quad (4)$$

Equation (2) shows that the firm demands are not modified between two successive planning steps. New firm demands (as they result from the consolidation process "g") remain consistent with their previous "flexible" values (3). The flexible bounds do not change between two planning steps (4).

### 4.2.2 Model of the supplier's behaviour

The *supplier's demand management process* computes the specification of the demand that will be taken into account in the supplier's planning process in its deterministic form. This management process depends on the uncertainty associated to the demand corresponding to the flexibility interval associated to the customer's demand. Thus, regarding the considered horizon (i.e. firm or flexible), the supplier will satisfy either equation 5 or 6.

$$\hat{D}_{p,t}^{\tau} = D_{p,t}^{\tau} \quad \forall p \quad \forall t \in FH^{\tau} \quad (5)$$

$$\hat{D}_{p,t}^{\tau} = f(\underline{D}_{p,t}^{\tau}, \overline{D}_{p,t}^{\tau}) \quad \forall p \quad \forall t \in LH^{\tau} \quad (6)$$

in which $\hat{D}_{p,t}^{\tau}$ is the deterministic demand on which the planning process is based. The definition of a value $\hat{D}_{p,t}^{\tau}$ is made using the demand management strategy $f$ as described in equation 6.

The planning behaviour is modelled as a planning problem using a mixed integer linear planning model (similar to those used in Advanced Planning Systems (APS)). Such a model is based on the one detailed in (Galasso *et al.*, 2006). The objective function (7) of this model has been adapted in order to maximise the gain calculated at each planning step. This model includes the following characteristics: multi-product, multi-components, possibility to adjust internal capacity through the use of extra-hours, change the workforce from one to two or three-shifts-work and subcontracting a part of the load.

The decision variables are introduced:
- $X_{p,t}$: internal production of final product $p$ at period $t$.
- $ST_{p,t}$: subcontracted production of final product $p$ at period $t$.
- $HS_t$: extra-hours used at period $t$.
- $B_{a,t}$: (binary variables) = 1 if action $a$ is used in order to modify the workforce at period $t$ and = 0 otherwise.

These decisions are linked with the following state variables:
- $I^+_{p,t}$ ; $I^-_{p,t}$: inventories and backorders levels at the end of period $t$ for the final product $p$.
- $J_{c,t}$: component inventory at period t.
- $A_{s,c,t}$: purchases of component $c$ bought at supplier $s$ to be delivered at period $t$.

The model is based on the following data:

- $CN$: nominal capacity available at each period $t$.
- $\{a\}$: set of actions that can be activated in order to adjust the capacity levels (i.e. 2 or 3-shifts-work) through the use of $B_{a,t}$.
- $\hat{D}_{p,t}$: deterministic demand of final product $p$ at period $t$ defined by the supplier.
- $R_p$: unitary production lead time for final product $p$.
- $\alpha_{p,c}$: bills of material coefficient linking final products $p$ and components $c$.

The planning model is defined hereafter in (7):

$$\max \sum_{t=\tau}^{\tau+HP-1} \left[ \sum_p v_p V_p - \sum_p h_p I^+_{p,t} - \sum_c c_c J_{c,t} - \sum_p b_p I^-_{p,t} - \sum_p u_p X_{p,t} - \sum_p st_p ST_{p,t} - \sum_c \sum_s f_{s,c} A_{s,c,t} - \sum_a o_a B_{a,t} - eHS_t \right] \quad (7)$$

is subject to:

$$I^+_{p,t} - I^-_{p,t} = I^+_{p,t-1} - I^-_{p,t-1} + X_{p,t-LP} + ST_{p,t-LS} - \hat{D}_{p,t} \quad \forall p, t \in HP \quad (8)$$

$$\sum_p R_p X_{p,t} \le CN + \sum_a (B_{a,t} \times SC_a) + HS_t \quad \forall t \in HP \quad (9)$$

$$J_{c,t} = J_{c,t-1} - \sum_p \alpha_{p,c}(X_{p,t} + ST_{p,t}) + \sum_s A_{s,c,t} \quad \forall t \in HP \quad (10)$$

$$\sum_p \alpha_{p,c}(X_{p,t} + ST_{p,t}) \le J_{c,t-1} \quad \forall t \in HP \quad (11)$$

$$HS_t \le HSMax \quad \forall t \in HP \quad (12)$$

The objective function (7) maximises the gain obtained through the plan established at each planning step. $v_p$ is the gain resulting from the deliveries of each product $p$. $h_p$, $hc_c$, $b_p$, $u_p$, $st_p$, $f_{s,c}$, $o_a$, $e$, are the unitary costs associated to the relevant decisions. Equation (8) links production quantities (subcontracted or not) and the levels of inventories and backorders. The lead times (LP standing for internal production and LS for the subcontracted production) are also introduced in equation (8). Moreover in that equation, the deterministic demand $\hat{D}_{p,t}$ is taken into account for the simulation of the planning process. The amount of production available at each period is limited by the capacity defined with constraint (9). A standard capacity $CN$ is available at each period. An amount of extra-hours $HS_t$ and an overcapacity $SC_a$ (introduced through the use of actions defined in $\{a\}$ and activated through $B_{a,t}$) can be added to the standard capacity. This constraint shows that resources are shared among products. Equation (10) enables the calculation of the inventory levels of components according to the purchases $A_{s,c,t}$ and the consumption of components linked to the internal and subcontracted production with the coefficients of the bills of materials $\alpha_{p,c}$. Constraint (11) ensure the consistency between the requirements and the inventory levels of components. Extra-hours are limited by (12) by a maximum value $HSMax$. All these constraints are defined at each period (time bucket) of the planning horizons. Each decision variable has its own dynamics and, similarly to the management of the customer demand, can be subject to a specific anticipation delay corresponding to the necessary organisational requirements previous to the applicability of such decisions.

## 5 ILLUSTRATIVE EXAMPLE

In this example, the cooperative decision making process detailed in section 4 is illustrated on an academic example (typical of aeronautics).

### 5.1 Decision and cooperation support approach (first run)

A first run of the specified decision and cooperation support approach is performed on the illustrative example.

#### 5.1.1 Problem and system definition

The system is a relationship with a customer and a SME supplier. In this system, the customer is supposed to order a single product representative of the aggregation at the tactical level of a family of end-items from the supplier. This product p is made of 1 component of type C1 ($\alpha_{p,C1}$ = 1) and of 2 components of type C2 ($\alpha_{p,C2}$ = 2). The supplier is in charge of assembling the two components in order to satisfy the demand specified and communicated by the customer.

The supply chain has been defined so that the length of the horizon on which the customer's demand is given enables the supplier to use all his decisional levers (i.e. use of extra-hours, subcontracting and use of both suppliers). This length encompasses the 4 periods necessary for the use of the subcontractor plus the four periods necessary to the use of the supplier 1 at rank 2 plus the 2 periods of the planning periodicity that is 12 periods.

The delays synthesised in Table 2 show several reactivity levels of each decision variables used in this behaviour model. The end-items internally produced can be delivered or added to the inventory at the following period ($LP$ = 1). The use of the subcontractor requires the transmission of the information 2 periods in advance ($DA_{ST}$ = 2), which forces the decisions to be frozen in the first two periods of the planning horizon. Then, the lead-time for the subcontractor is 2 periods ($LS$ = 2). Extra-hours must be anticipated with a delay of 1 period ($DA_{HS}$ = 1) and, obviously, are applied immediately. In order to make sure that rank 2 suppliers are able to manage their own production process, an anticipation of 4 periods is required for the rank 2 supplier 1 ($s_1$). An anticipation of 2 periods is required for the rank 2 supplier 2 ($s_2$). Thus, $s_2$ is more reactive than $s_1$.

| Decision | Lead time (LP) | Anticipation delay (DA) |
|---|---|---|
| Internal Production | 1 | |
| Subcontracting | 2 | 2 |
| Extra-hours | | 1 |
| Rank 2 Supplier 1 | | 4 |
| Rank 2 Supplier 2 | | 2 |

Table 2. Temporal parameters values

The unitary costs associated to each decision variable are given in Table 3:

| Decision variable | Unitary cost | Decision variable | Unitary cost |
|---|---|---|---|
| Purchase of C1 at S1 ($f_{S1,C1}$) | 2 | Backorders ($b_p$) | 20 |
| Purchase of C2 à S1 ($f_{S1,C2}$) | 1 | Final product holding ($h_p$) | 10 |
| Purchase of C1 à S2 ($f_{S2,C1}$) | 3 | Production ($u_p$) | 5 |
| Purchase of C2 à S2 ($f_{S2,C2}$) | 2 | Subcontracting ($st_p$) | 70 |
| Holding cost of C1 ($c_{C1}$) | 1 | Extra-hours ($e$) | 30 |
| Holding cost of C2 ($c_{C2}$) | 0,5 | | |

Table 3. Cost structure for the simulation

The selling price of final products p is 200cu and finally, the internal production of final product p requires 2 time units ($R_p$ = 2). It is interesting to notice that the supplier will need to choose among its suppliers in order to balance the need for the most reactive supplier (i.e. the rank 2 supplier 2) and minimising the purchasing cost as the first supplier is less expensive. Regarding the parameters associated to the production decision variables, the cost parameters privilege the use of internal production then, the use of extra-hours and finally, the use of the subcontractor.

In this context, the main problem of the *supplier* is to choose a planning strategy concerning the demand management in order to take into account the demand communicated by the customer in its planning process. Regarding the *customer* supply management process, an important decisional lever is the length of the firm and flexible horizon. Through this lever, the supplier has more or less visibility on the demand and thus more or less time to react and adapt its production process.

The objective of the supplier is to maximise the global gain using the best planning strategy according to the characteristics of its production process. On the other hand, the objective of the customer is to minimise the backorder levels. Both of them aim at agreeing on a common set of strategies dealing with these two goals while adjusting the model they use in order to simulate their planning and demand management process. Nevertheless, as the problem cannot be totally defined at once, other indicators have to be computed (stock levels, service level, production cost,…). At the customer level, the cost of backorders is considered as relevant.

### 5.1.2 Design of experiment

In our example, the customer identifies two possible trends for the evolution of his requirements. The identification of these trends shows the will of the customer to facilitate the organisation of its supplier. The uncertainty remaining on the customer requirement is characterised either by the possibility of occurrence of both trends and, moreover, by a flexibility of +/- 20% required for each trend.

The first trend (T1) reflects a strong punctual increase of the requirements with the acceptability of orders beyond the standard production capacity. Figure 5 shows the corresponding forecasts. The second trend (T2) presented in Figure 6 corresponds to a moderate increase as viewed by the customer. This punctual increase, expected for periods 20 to 25 is much lower than the previous one.

Both in figures 5 and 6, the minimum, the maximum and the average values of the demand are given at each period.

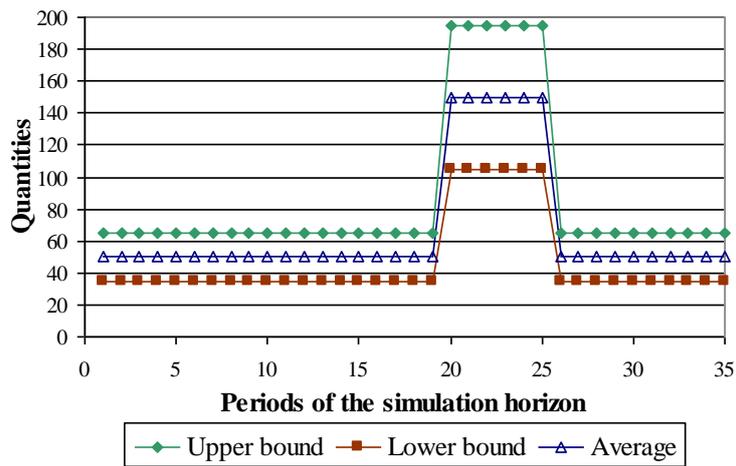
Figure 5. Trend 1

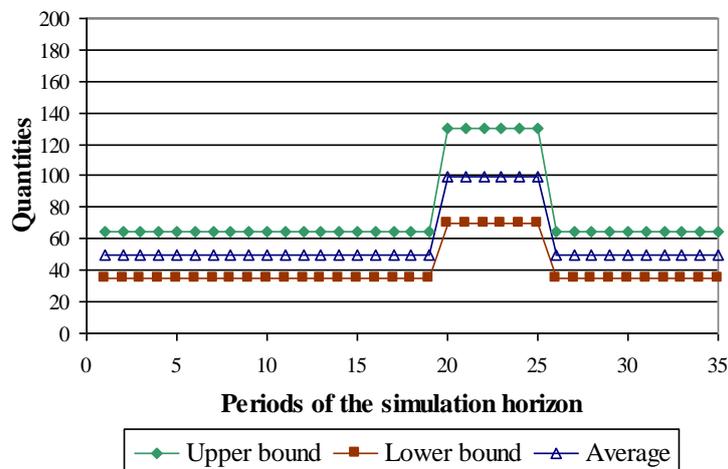
Figure 6. Trend 2

According to its height, the peak will have more or less influence on the planning process of the supplier and may require different uses of production capacities (internal with or without extra-hours, subcontracted) while taking into account the production delays.

In order to simulate several cooperative behavioural aspects, the demand management process of the customer includes a consolidation process for the demand and a visibility.

The consolidation process for the demand (noted g in Section 4.2.1) consists of an overestimation (resp. underestimation) of the demand noted "Min" (resp. "Max"). In that case, the customer will finally order the lower (resp. the upper) bound of the flexible demand.

Thus, the behaviour of the customer concerns the lengths of the firm horizon on which the firm demand is communicated by the customer to the supplier. Four lengths are studied so called visibilities:

- Visibility V1, on which the firm horizon length is 4 periods
- Visibility V2, on which the firm horizon length is 6 periods
- Visibility V3, on which the firm horizon length is 8 periods
- Visibility V4, on which the firm horizon length is 10 periods

Each of these lengths is completed by a flexible horizon the length of which constitutes the complement to the 12 periods of the planning horizon.

On the supplier side, as the planning process is run, the understanding of the trend increases. In order to manage the uncertainty on the flexible demand, the supplier uses two planning strategies (noted $f$ in section 4.2.1) S1 and S2, in its demand management process:

- S1: choose the maximum of the flexible demand $f = \max\left[\underline{D}_{p,t}^{\tau}, \overline{D}_{p,t}^{\tau}\right]$

- S2: choose the minimum of the flexible demand $f = \min\left[\underline{D}_{p,t}^{\tau}, \overline{D}_{p,t}^{\tau}\right]$

The evaluation of these strategies according to the trends and the consolidation processes is done running simulations that are designed as a combination of:

- a trend of the evolution of the customer's requirements (T1 or T2),
- a type of demand management for the customer :
  - behaviour "$g$" for the customer (overestimation denoted "Min" or under-estimation denoted "Max" of the demand),
  - visibility taken from the four lengths of the firm horizon communicated by the customer.
- a planning strategy "$f$" of the supplier (concerning the choice of the maximal flexible demand denoted S1 or the choice of the minimal denoted S2).

Thus, our design of experiments consists of 32 experiments in which cost and temporal parameters remain constant for each simulation.

In the next sections, the results for the supplier and the customer risk evaluation are detailed.

### 5.1.3 Simulations and computation of the indicators

The gains and the costs obtained during the simulations with the use of the strategy S1 (i.e. the supplier integrates the maximum values of the demand) and S2 (i.e. the supplier integrates the minimum values of the demand) are presented in Table 4.

| Strategy | Trend | Consolidation Process | Visibility | Global Gains | Global costs | Total Production cost | Total Inventory cost | Total Backorder cost | Total Purchasing cost |
|---|---|---|---|---|---|---|---|---|---|
| S1 | T1 | Min | V1 | 245 201 | 83 800 | 27 580 | 48 303 | 0 | 7 917 |
| | | | V2 | 275 477 | 53 523 | 19 585 | 25 976 | 300 | 7 662 |
| | | | V3 | 287 509 | 41 491 | 19 195 | 13 995 | 900 | 7 401 |
| | | | V4 | 291 328 | 37 673 | 20 970 | 8 408 | 900 | 7 395 |
| | | Max | V1 | 476 378 | 134 622 | 93 160 | 12 581 | 14 260 | 14 621 |
| | | | V2 | 477 185 | 133 816 | 93 070 | 12 573 | 13 620 | 14 553 |
| | | | V3 | 478 565 | 132 436 | 93 010 | 12 573 | 12 300 | 14 553 |
| | | | V4 | 478 610 | 132 391 | 93 010 | 12 573 | 12 300 | 14 508 |
| | T2 | Min | V1 | 235 470 | 51 530 | 9 840 | 34 940 | 0 | 6 750 |
| | | | V2 | 256 284 | 30 716 | 6 850 | 17 116 | 300 | 6 450 |
| | | | V3 | 262 128 | 24 873 | 6 850 | 11 178 | 500 | 6 345 |
| | | | V4 | 264 557 | 22 443 | 9 250 | 6 360 | 500 | 6 333 |
| | | Max | V1 | 444 191 | 88 809 | 66 130 | 6 356 | 3 760 | 12 563 |
| | | | V2 | 446 378 | 86 623 | 65 980 | 6 348 | 1 800 | 12 495 |
| | | | V3 | 446 378 | 86 623 | 65 980 | 6 348 | 1 800 | 12 495 |
| | | | V4 | 446 423 | 86 578 | 65 980 | 6 348 | 1 800 | 12 450 |
| S2 | T1 | Min | V1 | 291 798 | 37 203 | 21 020 | 6 973 | 1 800 | 7 410 |
| | | | V2 | 291 798 | 37 203 | 21 020 | 6 973 | 1 800 | 7 410 |
| | | | V3 | 291 798 | 37 203 | 21 020 | 6 973 | 1 800 | 7 410 |
| | | | V4 | 291 798 | 37 203 | 21 020 | 6 973 | 1 800 | 7 410 |
| | | Max | V1 | 403 344 | 207 657 | 88 040 | 7 323 | 96 040 | 16 254 |
| | | | V2 | 444 947 | 166 054 | 90 685 | 7 413 | 52 140 | 15 816 |
| | | | V3 | 463 995 | 147 006 | 93 010 | 7 933 | 30 700 | 15 363 |
| | | | V4 | 473 611 | 137 390 | 93 010 | 9 803 | 19 940 | 14 637 |
| | T2 | Min | V1 | 264 853 | 22 148 | 9 550 | 5 198 | 1 100 | 6 300 |
| | | | V2 | 264 853 | 22 148 | 9 550 | 5 198 | 1 100 | 6 300 |
| | | | V3 | 264 853 | 22 148 | 9 550 | 5 198 | 1 100 | 6 300 |
| | | | V4 | 264 853 | 22 148 | 9 550 | 5 198 | 1 100 | 6 300 |
| | | Max | V1 | 383 765 | 149 236 | 60 895 | 6 348 | 67 500 | 14 493 |
| | | | V2 | 425 302 | 107 699 | 63 670 | 6 348 | 23 260 | 14 421 |
| | | | V3 | 444 929 | 88 072 | 65 935 | 6 348 | 2 100 | 13 689 |
| | | | V4 | 446 378 | 86 623 | 65 980 | 6 348 | 1 800 | 12 495 |

Table 4. Gains obtained for S1 and S2

| Visibility | Trend | Customer behaviour | Strategy | Total Backorder cost | Global Gains | Global costs | Total Production cost | Total Inventory cost | Total Purchasing cost |
|---|---|---|---|---|---|---|---|---|---|
| V1 | T1 | min | S1 | 0 | 245 201 | 83 800 | 27 580 | 48 303 | 7 917 |
| | | | S2 | 1 800 | 291 798 | 37 203 | 21 020 | 6 973 | 7 410 |
| | | Max | S1 | 14 260 | 476 378 | 134 622 | 93 160 | 12 581 | 14 621 |
| | | | S2 | 96 040 | 403 344 | 207 657 | 88 040 | 7 323 | 16 254 |
| | T2 | min | S1 | 0 | 235 470 | 51 530 | 9 840 | 34 940 | 6 750 |
| | | | S2 | 1 100 | 264 853 | 22 148 | 9 550 | 5 198 | 6 300 |
| | | Max | S1 | 3 760 | 444 191 | 88 809 | 66 130 | 6 356 | 12 563 |
| | | | S2 | 67 500 | 383 765 | 149 236 | 60 895 | 6 348 | 14 493 |
| V2 | T1 | min | S1 | 300 | 275 477 | 53 523 | 19 585 | 25 976 | 7 662 |
| | | | S2 | 1 800 | 291 798 | 37 203 | 21 020 | 6 973 | 7 410 |
| | | Max | S1 | 13 620 | 477 185 | 133 816 | 93 070 | 12 573 | 14 553 |
| | | | S2 | 52 140 | 444 947 | 166 054 | 90 685 | 7 413 | 15 816 |
| | T2 | min | S1 | 300 | 256 284 | 30 716 | 6 850 | 17 116 | 6 450 |
| | | | S2 | 1 100 | 264 853 | 22 148 | 9 550 | 5 198 | 6 300 |
| | | Max | S1 | 1 800 | 446 378 | 86 623 | 65 980 | 6 348 | 12 495 |
| | | | S2 | 23 260 | 425 302 | 107 699 | 63 670 | 6 348 | 14 421 |
| V3 | T1 | min | S1 | 900 | 287 509 | 41 491 | 19 195 | 13 995 | 7 401 |
| | | | S2 | 1 800 | 291 798 | 37 203 | 21 020 | 6 973 | 7 410 |
| | | Max | S1 | 12 300 | 478 565 | 132 436 | 93 010 | 12 573 | 14 553 |
| | | | S2 | 30 700 | 463 995 | 147 006 | 93 010 | 7 933 | 15 363 |
| | T2 | min | S1 | 500 | 262 128 | 24 873 | 6 850 | 11 178 | 6 345 |
| | | | S2 | 1 100 | 264 853 | 22 148 | 9 550 | 5 198 | 6 300 |
| | | Max | S1 | 1 800 | 446 378 | 86 623 | 65 980 | 6 348 | 12 495 |
| | | | S2 | 2 100 | 444 929 | 88 072 | 65 935 | 6 348 | 13 689 |
| V4 | T1 | min | S1 | 900 | 291 328 | 37 673 | 20 970 | 8 408 | 7 395 |
| | | | S2 | 1 800 | 291 798 | 37 203 | 21 020 | 6 973 | 7 410 |
| | | Max | S1 | 12 300 | 478 610 | 132 391 | 93 010 | 12 573 | 14 508 |
| | | | S2 | 19 940 | 473 611 | 137 390 | 93 010 | 9 803 | 14 637 |
| | T2 | min | S1 | 500 | 264 557 | 22 443 | 9 250 | 6 360 | 6 333 |
| | | | S2 | 1 100 | 264 853 | 22 148 | 9 550 | 5 198 | 6 300 |
| | | Max | S1 | 1 800 | 446 423 | 86 578 | 65 980 | 6 348 | 12 450 |
| | | | S2 | 1 800 | 446 378 | 86 623 | 65 980 | 6 348 | 12 495 |

Table 5. Backorder costs obtained for V1, V2, V3 and V4

### 5.1.4 Risk evaluation

The risk evaluation approach step is run from the supplier's point of view then and the customer's point of view.

#### 5.1.4.1 Risk evaluation from the supplier's point of view

The best gains obtained for each behaviour (S1 or S2) of the supplier are (cf. Table 4): 478 610 for the strategy S1 and of 473 611 for the strategy S2. The worst gains are 235 470 for S1 and 264 853 for S2. According to these results, it is possible to derive the risk diagram presented in Figure 8. To do so, it is necessary to calculate the value of the coefficient of optimism α of the Hurwicz criterion beyond which a change of strategy is "recommended" (cf. Figure 7).

In order to visualise this specific point, we plot the straight lines corresponding to equations:
- $H_{S1} = (1-\alpha)\times 235\,470 + \alpha\times 478\,610$ for S1 and
- $H_{S2} = (1-\alpha)\times 264\,853 + \alpha\times 403\,344$ for S2.

From these equations, it is easily possible to determine that $H_{S1} = H_{S2}$ for $\alpha_1 \approx 0.855$. It is now possible to establish the risk diagram (Figure 8). Firstly the α-axis symbolising the propensity to risk of the decision maker is drawn, highlighting the value of parameter α indicating a change of strategy (here for α = 0,855). Then, both strategies S1 and S2 are placed on the axis. Finally, the other criteria (Laplace and Savage) are superposed to the strategy that they recommend in the diagram.

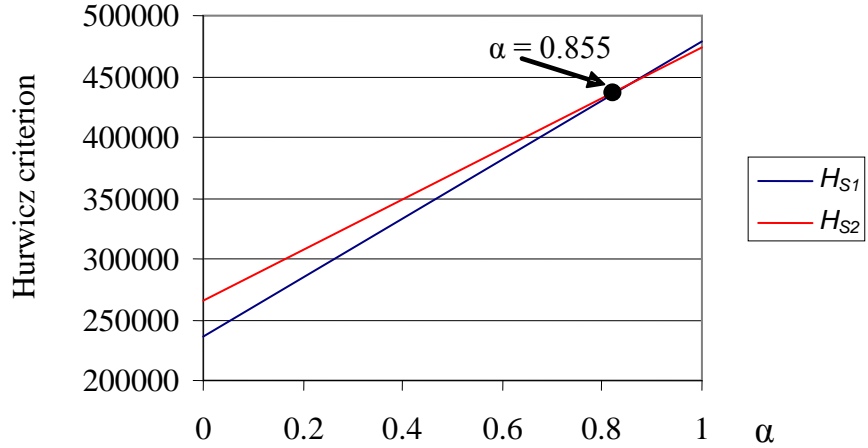

Figure 7. Point of change of supplier's strategy

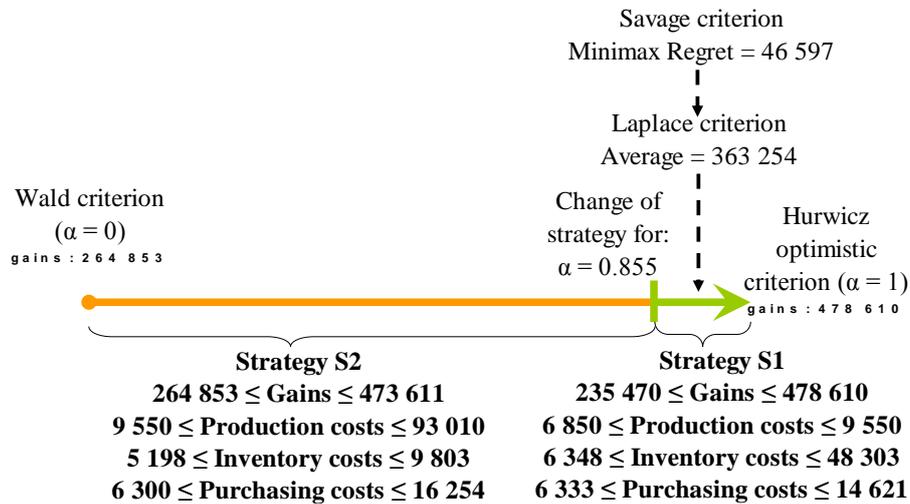

Figure 8. Risk diagram for the supplier's point of view

In order to facilitate the decision making, Table 6 and figures 9 to 12 are added.

|  | Regret using S1 | Regret using S2 |
| --- | --- | --- |
| Min regret / S1 | 0 | -46 597 |
| Max regret / S1 | 0 | 73 034 |
| Min regret / S2 | -73 034 | 0 |
| Max regret / S2 | 46 597 | 0 |

Table 6. Regrets table for S1 and S2

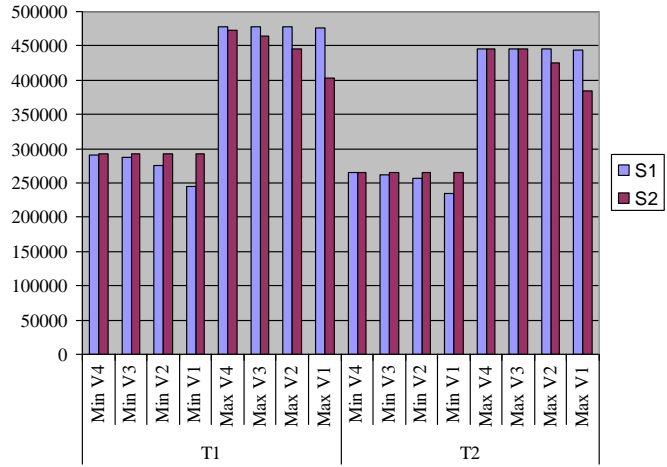

Figure 9. Evolution of gains

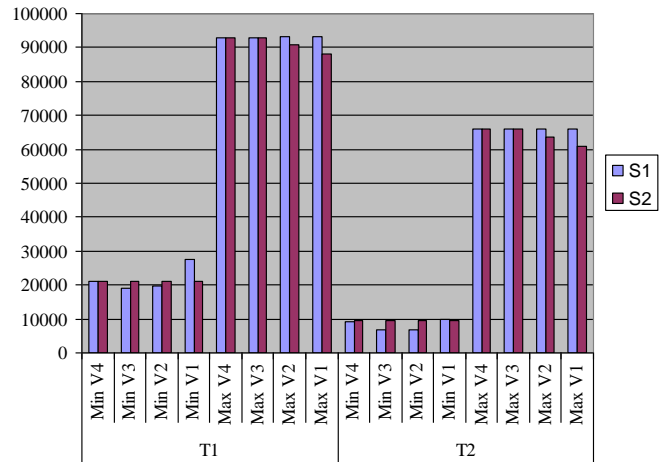

Figure 10. Evolution of production costs

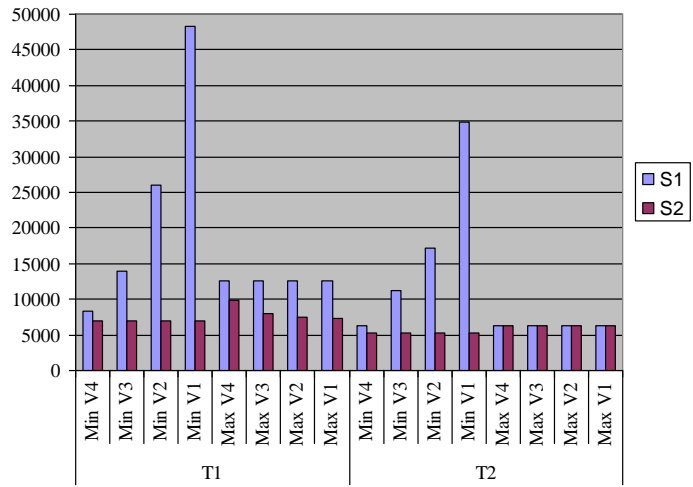

Figure 11. Evolution of inventory costs

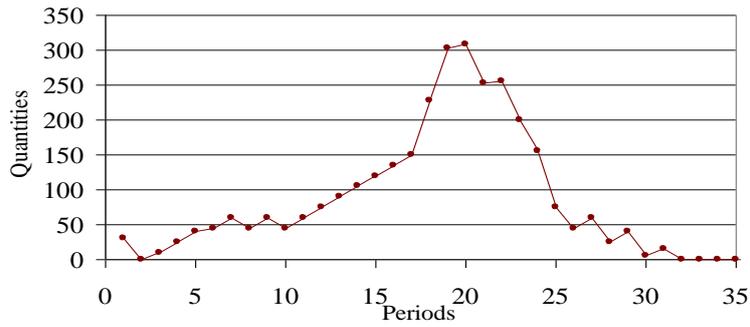

Figure 12. Evolution of inventories in the highest inventory costs case

#### 5.1.4.2 Risk evaluation from the customer's point of view

In the same way, Figure 13 shows, from the customer point of view, the values of the Hurwicz criterion for the four visibilities that could be given by the customer. This figure is magnified in a range of values of α comprised between 0.95 and 1 in order to make the intersection point more visible.

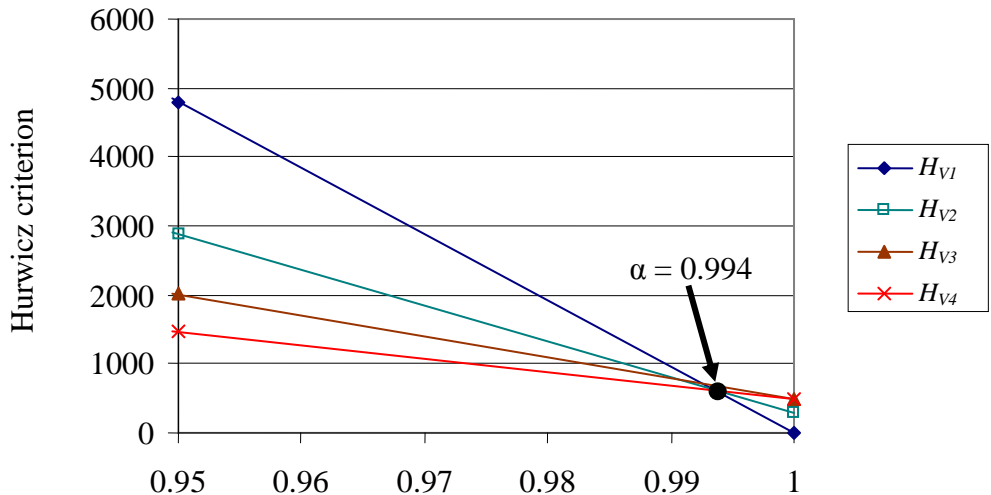

Figure 13. Point of change of customer's strategy

This diagram in Figure 13 permits to establish the following risk diagram in Figure 14.

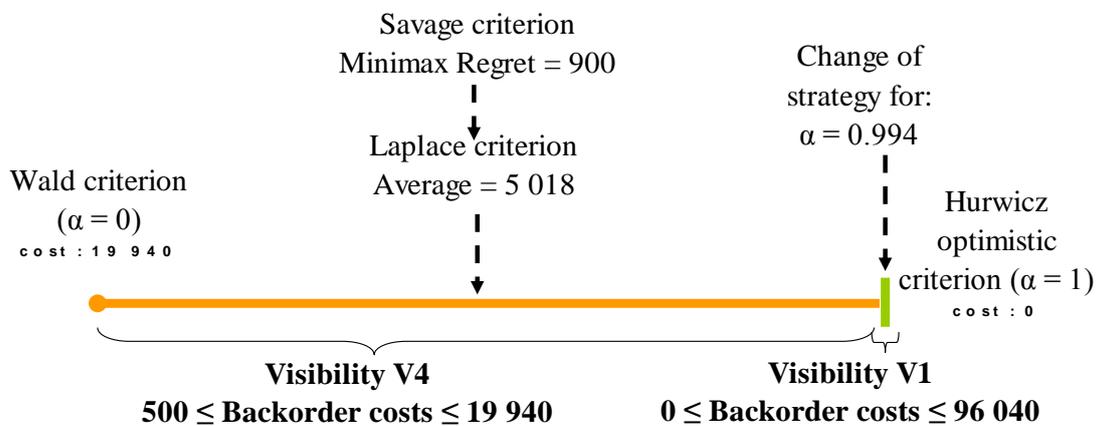

Figure 14. Risk diagram from the customer's point of view

In order to facilitate the decision making, Table 7 and figures 15 and 16 are added.

|                  | Regret using V1 | Regret using V2 | Regret using V3 | Regret using V4 |
|------------------|-----------------|-----------------|-----------------|-----------------|
| Min regret / V1  | 0               | -44 240         | -65 400         | -76 100         |
| Max regret / V1  | 0               | 300             | 900             | 900             |
| Min regret / V2  | -300            | 0               | -21 440         | -32 200         |
| Max regret / V2  | 44 240          | 0               | 600             | 600             |
| Min regret / V3  | -900            | -600            | 0               | -10 760         |
| Max regret / V3  | 65 400          | 21 440          | 0               | 0               |
| Min regret / V4  | -900            | -600            | 0               | 0               |
| Max regret / V4  | 76 100          | 32 200          | 10 760          | 0               |

Table 7. Regrets using V1, V2, V3 or V4

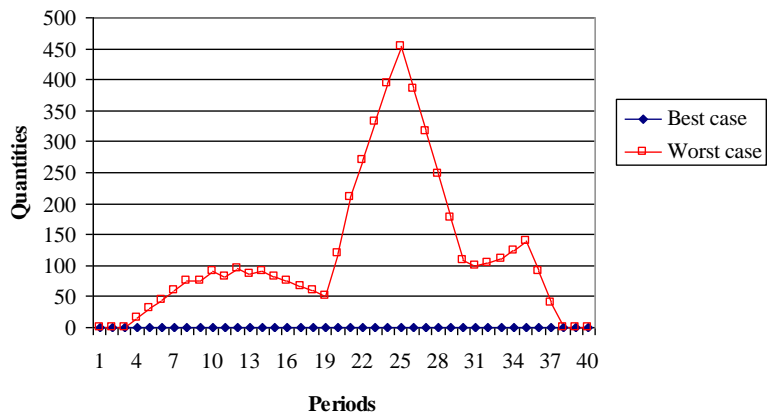

Figure 15. Backorders evolution for V1

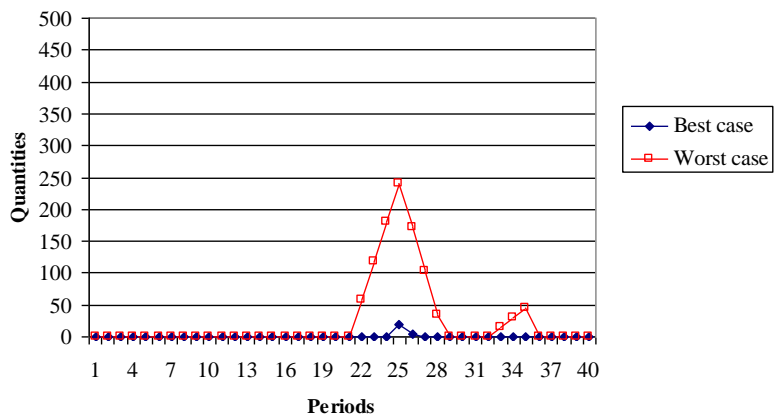

Figure 16. Backorders evolution for V4

### 5.1.5 Design of cooperative processes

After the evaluation step, a risk analysis is performed from the two actors' point of view.

#### 5.1.5.1 Risk analysis from the supplier's point of view

The risk diagram in Figure 8 shows that when a pessimistic point of view is adopted (α tends to 0) the planning strategy using the minimal demand (S2) is recommended. The weighted Hurwicz criterion proposes a change in the strategy applied for an optimism degree of 0.855 (values comprised between 0 and 1). This value means that the strategy S2 may be envisaged by the supplier's decision maker even if other criteria such as Laplace or Savage recommend the choice of the strategy S1. S1 is also recommended by the Hurwicz criteria for values over α = 0.855. Thus, it is the interest of the supplier's decision maker to request additional information (i.e. information from the customer or concerning the global market evolution) in order to determine if he should be pessimistic or not. For each recommended strategy, the minimal and maximal values of the production, inventory and purchasing costs are indicated. These results provide further appeal to a simple simulation giving raw gains according to several scenarios. Indeed, in a first approach, it could be obvious that the higher the demand is, the higher the gains are. Nevertheless, disruptions may put into question the occurrence of a scenario leading to such gains and the raw results remain uncertain. Therefore, through the risk diagram, we provide not only information regarding an interesting strategy to be applied but also an indication about the relevance of this choice.

Table 6 indicates clearly that the regret of using S2 instead of S1 is the most important. This statement seems contradictory regarding the risk diagram in Figure 8 which recommends the strategy S2. Thus, it is important to corroborate this result with the other indicators (figures 9 to 12).

From Figure 11, it can be pointed out that the strategy S1 can generate high inventory costs (i.e. 48 303cu.). Considering these costs, the supplier can found an incentive to modify its Advanced Planning System in order to integrate a maximal level of inventory, which corresponds to the maximum costs he would like to allocate to inventories.

### 5.1.5.2   Risk analysis from the customer's point of view

It can be extracted from the risk diagram in Figure 14 that the visibility V4 is the best one for a large part in Hurwicz sense. Thus it can be privileged for values of α besides 0.994. Beyond this value, visibility V1 can be privileged for a small range that can hardly be seen in the diagram.

Considering the regret table associated to the risk diagram and given in Table 7, it can be stressed that the maximum regret of each strategy is obtained in the case of a substitution with visibility V4. Without considering the visibility V4, the maximum regret of using V3 instead of visibilities V1 and V2 is limited to 900cu while V1 generates a maximal regret of 65 400cu and V2 of 21 400cu.

This information is confirmed by figures 15 and 16 provided in the dashboard and detailing the evolution of backorders for the best and worst cases for each visibility. Thus, the impact of the visibility over the backorders is clearly identified and demonstrated. The peak of backorders at period 25 is reduced when more visibility is given to the supplier.

### 5.1.5.3   Risk analysis from the cooperation point of view

Finally, the dashboard leads the decision makers to determine the pair of strategies S2 for the supplier and V4 for the customer. This evaluation assumes that the addition of visibility is free of cost for the customer. But increasing visibility requires an effort for the customer that has to assume firm requirements communicated to the supplier over a more important horizon. Thus, the customer has to consolidate his own planning and may be production processes in order to make such increase of visibility feasible. This necessary effort can be integrated in the evaluation process through penalty costs.

## 5.2   Decision and cooperation approach (second run)

In order to consider the modifications of the system that have been suggested thanks to the first run of the specified approach, a second run is performed.

### 5.2.1   Problem and system definition (second run)

Regarding the previous experiments in sections 5.3 and 5.4, two perspectives of improvement of the system have been pointed out. The first one comes from the supplier's viewpoint and is linked to the storage capacity that should be limited. The second one is linked to the customer's viewpoint and reveals a need for the integration of effort required in order to increase the visibility given to the supplier. This storage capacity and a measure of the effort of the customer evaluated in terms of costs are defined according to the decision maker knowledge. The adjustment of these values through "what-if" simulations is part of the improvement process of the decision and cooperation support approach as well as the integration of these recommendations.

In our case, the supplier may consider in a second approach that a maximum level of 80 products could be investigated. This value is considered as a first improvement of the model and could be adjusted according to the results of a new global simulation run. Moreover, concerning the customer point of view, the necessary effort in order to increase the visibility has been made through the addition of penalty costs as given in Table 8.

|  | V1 to V2 | V1 to V3 | V1 to V4 |
|---|---|---|---|
| Penalty costs | 1000 | 2000 | 5000 |

Table 8. Penalty costs for increasing visibility

### 5.2.2 Design of experiment (second run)

After the redefinition of the problem under study, the same complete set of experiments is used.

### 5.2.3 Simulations and computation of the indicators (second run)

The results of these simulations are given in Table 9 below.

| Strategy | Trend | Consolidation Process | Visibility | Global Gains | Global costs | Total Production cost | Total Inventory cost | Total Backorder cost | Total Purchasing cost |
|---|---|---|---|---|---|---|---|---|---|
| S1 | T1 | Min | V1 | 259 783 | 69 218 | 30 065 | 31 353 | 0 | 7 800 |
|  |  |  | V2 | 281 771 | 47 229 | 23 130 | 16 299 | 1 300 | 7 500 |
|  |  |  | V3 | 288 998 | 40 003 | 21 675 | 10 033 | 2 900 | 7 395 |
|  |  |  | V4 | 291 373 | 37 628 | 21 600 | 7 733 | 5 900 | 7 395 |
|  |  | Max | V1 | 465 701 | 145 299 | 93 115 | 8 431 | 29 240 | 14 513 |
|  |  |  | V2 | 466 508 | 144 493 | 93 025 | 8 423 | 29 600 | 14 445 |
|  |  |  | V3 | 467 888 | 143 113 | 92 965 | 8 423 | 29 280 | 14 445 |
|  |  |  | V4 | 467 933 | 143 068 | 92 965 | 8 423 | 32 280 | 14 400 |
|  | T2 | Min | V1 | 236 485 | 50 515 | 14 130 | 29 635 | 0 | 6 750 |
|  |  |  | V2 | 256 946 | 30 054 | 8 350 | 14 954 | 1 300 | 6 450 |
|  |  |  | V3 | 262 815 | 24 185 | 8 050 | 9 290 | 2 500 | 6 345 |
|  |  |  | V4 | 264 557 | 22 443 | 9 250 | 6 360 | 5 500 | 6 333 |
|  |  | Max | V1 | 444 191 | 88 809 | 66 130 | 6 356 | 3 760 | 12 563 |
|  |  |  | V2 | 444 998 | 88 003 | 66 040 | 6 348 | 4 120 | 12 495 |
|  |  |  | V3 | 446 378 | 86 623 | 65 980 | 6 348 | 3 800 | 12 495 |
|  |  |  | V4 | 446 423 | 86 578 | 65 980 | 6 348 | 6 800 | 12 450 |
| S2 | T1 | Min | V1 | 291 798 | 37 203 | 21 020 | 6 973 | 1 800 | 7 410 |
|  |  |  | V2 | 291 798 | 37 203 | 21 020 | 6 973 | 2 800 | 7 410 |
|  |  |  | V3 | 291 798 | 37 203 | 21 020 | 6 973 | 3 800 | 7 410 |
|  |  |  | V4 | 291 798 | 37 203 | 21 020 | 6 973 | 6 800 | 7 410 |
|  |  | Max | V1 | 403 344 | 207 657 | 88 040 | 7 323 | 96 040 | 16 254 |
|  |  |  | V2 | 444 965 | 166 036 | 90 685 | 7 413 | 53 140 | 15 798 |
|  |  |  | V3 | 463 995 | 147 006 | 93 010 | 7 933 | 32 700 | 15 363 |
|  |  |  | V4 | 467 861 | 143 140 | 92 965 | 8 423 | 32 280 | 14 472 |
|  | T2 | Min | V1 | 264 853 | 22 148 | 9 550 | 5 198 | 1 100 | 6 300 |
|  |  |  | V2 | 264 853 | 22 148 | 9 550 | 5 198 | 2 100 | 6 300 |
|  |  |  | V3 | 264 853 | 22 148 | 9 550 | 5 198 | 3 100 | 6 300 |
|  |  |  | V4 | 264 853 | 22 148 | 9 550 | 5 198 | 6 100 | 6 300 |
|  |  | Max | V1 | 383 720 | 149 281 | 60 895 | 6 348 | 67 500 | 14 538 |
|  |  |  | V2 | 425 302 | 107 699 | 63 670 | 6 348 | 24 260 | 14 421 |
|  |  |  | V3 | 444 974 | 88 027 | 65 935 | 6 348 | 4 100 | 13 644 |
|  |  |  | V4 | 446 378 | 86 623 | 65 980 | 6 348 | 6 800 | 12 495 |

Table 9. Gains and costs obtained using limited inventories and penalty costs

### 5.2.4 Risk evaluation (second run)

A new dashboard is provided to the decision makers (Figures 17 to 24 and Table 10).

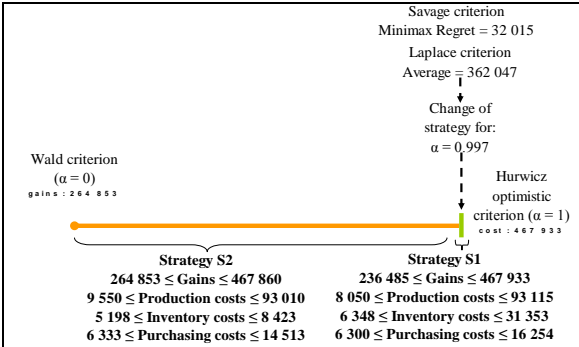

Figure 17. Risk diagram for the supplier's point of view using limited inventories

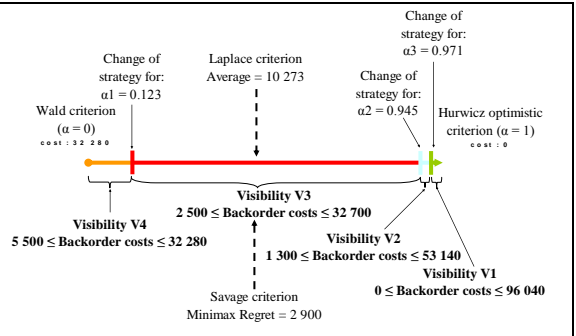

Figure 21. Risk diagram from the customer's point of view with penalties

|  | Regret using V1 | Regret using V2 | Regret using V3 | Regret using V4 |
|---|---|---|---|---|
| Min regret / V1 | 0 | -43 240 | -63 400 | -63 760 |
| Max regret / V1 | 0 | 1 300 | 2 900 | 5 900 |
| Min regret / V2 | -1 300 | 0 | -20 440 | -20 860 |
| Max regret / V2 | 43 240 | 0 | 1 600 | 4 600 |
| Min regret / V3 | -2 900 | -1 600 | 0 | -420 |
| Max regret / V3 | 63 400 | 20 440 | 0 | 3 000 |
| Min regret / V4 | -5 900 | -4 600 | -3 000 | 0 |
| Max regret / V4 | 63 760 | 20 860 | 420 | 0 |

Table 10. Regrets using V1, V2, V3 or V4 with penalties

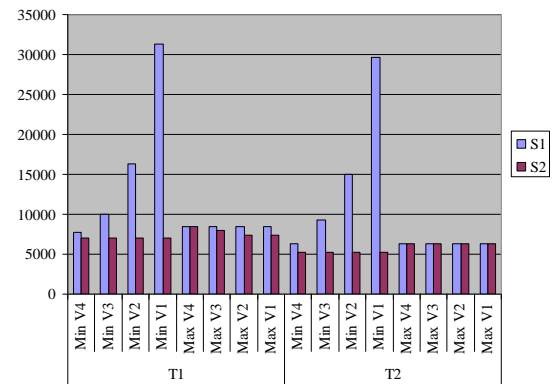

Figure 18. Evolution of inventory costs using limited inventories

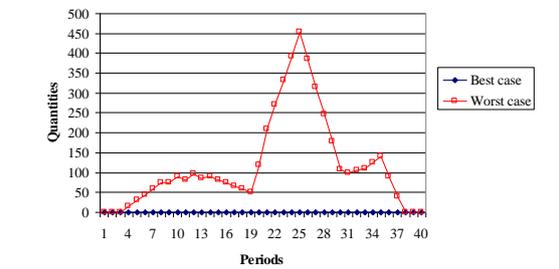

Figure 22. Backorders evolution for V1

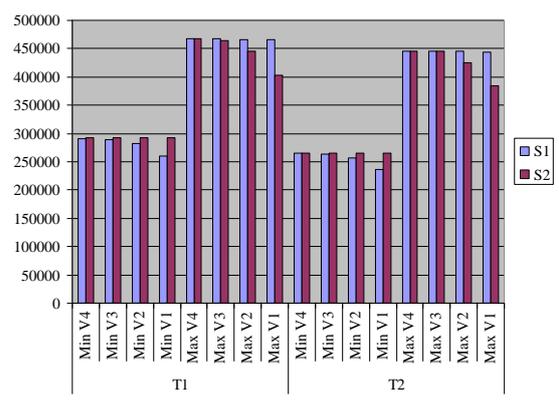

Figure 19. Evolution of gains using limited inventories

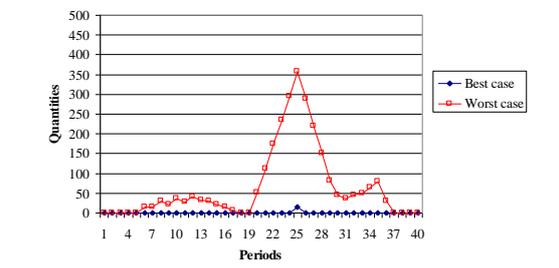

Figure 23. Backorders evolution for V2

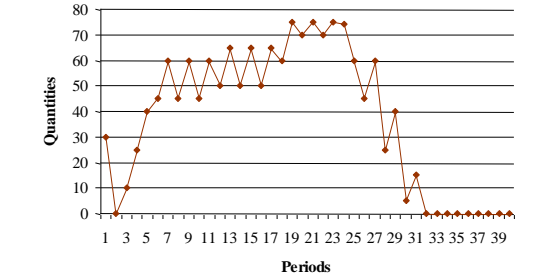

Figure 20. Evolution of limited inventories in the highest inventory costs case

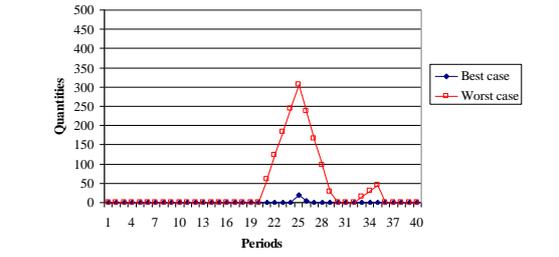

Figure 24. Backorders evolution for V3 and V4

### 5.2.5 Design of cooperative processes (second run)

After the evaluation step, a risk analysis (second run) is performed, which leads to the design of cooperative processes.

#### 5.2.5.1 Risk analysis from the supplier's point of view (second run)

Considering Figure 17, the risk diagram in Figure 8 is significantly modified and the use of the strategy S2 for the supplier is more recommended. Figure 18 shows that the evolution of inventory costs highlights the use of the stock limitation as the highest inventory costs significantly decrease.
Figure 19 shows the evolution of the gains using limited inventory. It stresses that the gap between the two strategies S1 and S2 slightly decreases. It is not surprising since the stock level has been constrained. Hence the anticipation effect provided by S1 is reduced as well.

#### 5.2.5.2 Risk analysis from the customer's point of view (second run)

The modifications given to the evaluation of risk has altered the view of the customer. Considering the risk diagram in Figure 14, visibility V3 now appears for an important range of values of α and will probably be chosen by the customer.

#### 5.2.5.3 Risk analysis from the cooperation's point of view (second run)

While considering this new definition of the production problem and the links between the strategies available for the customer, the two dashboards lead the decision makers to determine the pair of strategies S2 for the supplier and V3 for the customer.

#### 5.2.5.4 Conclusion of the partnership contract

It is important to notice that the determination of this pair (S2 , V3) requires the cooperation of the two decision makers in order to provide reliable data, share the analysis of results, improve the modelling process and finally, apply the strategies on both sides.

## 6 CONCLUSION

In this article, a decision and cooperation support approach for the design of cooperative planning process is proposed which is dedicated to a customer-supplier relationship.
In this approach, after a design step for the problem under study, the parameters that arise from this step are instantiated in order to generate the experiments to be simulated.
Then the exploitation of these simulations is done using a dedicated dashboard. This dashboard includes risks diagrams built according to the weighted Hurwicz criterion and other criteria (i.e. Laplace, Wald and Savage). These diagrams give more information than a simple evaluation of the plans established by the supplier according to the demand given by the customer. Indeed, thanks to the Hurwicz criterion, they introduce degrees of optimism for which a planning strategy can be privileged for each actor. Moreover other indicators are provided to decision makers including a regret table which situates the strategies proposed by the risk diagram within the set of potential strategies and other evolution curves (inventory, production, purchasing…)
In that way, decision makers can cooperatively define the planning strategies at the supplier's level and demand management strategies at the customer level.
The ability of our approach to evaluate and compare the impact of these strategies simultaneously provides decision makers with tools for using their expert knowledge.

Finally, the illustrative example shows the importance of viewing these simulations as a cooperative process: adjustments of both models and parameters values may be necessary after a first run of experiments. These adjustments require the experience of decisions makers from both the supplier and the customer sides.

There are many perspectives to this work. Firstly, more configurations of the planning parameters and of the demand that is communicated by the customer to the supplier should be investigated. We expect a confirmation of performance improvement when the cooperation level in the customer-supplier relationship is improved. Furthermore, an extension to linear or networked supply chains should be investigated. Thus, we may obtain a set of strategies that can be used at each rank of the chain while improving its global performance.

**REFERENCES**

Affonso, R., Marcotte, F., Grabot, B., 2006. Coordination model in the supply chain, *2nd I*PROMS virtual conference*, Pham, Eldukhri and Soroka (Eds.), Elsevier, pp. 468-473

Bartezzaghi, E., Verganti, R., 1995. Managing demand uncertainty through order overplanning. *International Journal of Production Economics*, 40, pp. 107-120.

Bouchon-Meunier, B., 1995. *La logique floue et ses applications* Addison-Weyley.

Bräutigam, J., Mehler-Bicher, A., Esche, C., 2003. Uncertainty as a key value driver of real options. *Proceedings of the annual international conference on real options*, available online at: http://www.realoptions.org/papers2003/BraeutigamUncertainty.pdf.

Brindley, C. (Ed.), 2004. *Supply chain risk*. MPG Books Ltd.

Christopher, M., 2003. Understanding Supply Chain Risk: A Self-Assessment Workbook. Cranfield University, School of Management, http://www.som.cranfield.ac.uk/som/research/centres/lscm/risk2002.asp , visited January 17, 2008.

Dubois, D., Prade, H., 1988. *Possibility Theory: An Approach to Computerized Processing of Uncertainty*, Plenum Press, New York.

Dudek, G., Stadtler, H., 2005. Negotiation-based collaborative planning between supply chains partners, *European Journal of Operational Research*, 163(3), pp. 668-687.

Galasso, F., Mercé, C., Grabot, B., 2006. Decision support for supply chain planning under uncertainty. *Proceedings from the 12$^{th}$ IFAC International Symposium Information Control Problems in Manufacturing (INCOM),* St-Etienne, France, 3, pp. 233-238.

Ganeshan, R., Jack, R., Magazine M.J., Stephens, P., 1999. A Taxonomic Review of Supply Chain Management Research, *in Quantitative Models for Supply Chain Management*, Kluwer Academic Publishers, Boston, pp. 841-880.

Génin, P., Thomas, A., Lamouri, S., 2007. How to manage robust tactical planning with an APS (Advanced Planning Systems). *Journal of Intelligent Manufacturing*, 18, pp. 209-221.

Holton, G.A., 2004. Defining Risk. *Financial Analysts Journal,* 60(6), pp. 19-25.

Huang, G.Q., Lau, J.S.K., Mak, K.L., 2003. The impacts of sharing production information on supply chain dynamics: a review of the literature. *International Journal of Production Research*, 41(7), pp. 1483-1517.

Lamothe, J., Mahmoudi, J., Thierry, C., 2008. *Cooperation to reduce risk in a telecom supply chain.* in: Supply Chain Forum, IADIS Digital Library, Special edition of Managing Supply Chain Risk, (in press).

Lang, J., 2003. *Contribution à l'étude de modèles, de langages et d'algorithmes pour le raisonnement et la prise de décision en intelligence artificielle.* Habilitation à Diriger des Recherches, Université Paul Sabatier, France.

Lapide, L., 2001. New developments in business forecasting, *Journal of Business Forecasting Methods & Systems*, 20(4), pp. 11-13.

Lee, H.L., Padmanabhan, P., Whang. S., 1997. Information Distortion in a Supply Chain: The Bullwhip Effect, *Management Science*, 43, pp. 546-558.

Mahmoudi, J., 2006. *Simulation et gestion des risques en planification distribuée de chaînes logistiques : application au secteur de l'électronique et des télécommunications.* Thèse de Doctorat, Sup'aero, France, In french.

McCarthy, T., Golicic, S., 2002. Implementing collaborative forecasting to improve supply chain performance, *International Journal of Physical Distribution & Logistics Management*, 32(6), pp. 431-454.

Moyaux, T., 2004. *Design, simulation and analysis of collaborative strategies in multi-agent systems: The case of supply chain management*, PhD thesis, Université Laval, Ville de Québec, Québec, Canada.

Mula, J., Poler, R., Garcia, J.P., Lario, F.C., 2006. Models for production planning under uncertainty: A review, *International Journal of Production Economics,* 103(1), pp. 271-285.

Odette, 2007. http://www.odette.org, visited January 17, 2008.

Ritchie, B., Brindley, C., 2004. *Risk characteristics of the supply chain – a contingency framework*. Brindley, C. (Ed.). In Supply chain risk. Cornwall: MPG Books Ltd.

Rota, K., Thierry, C., Bel, G., 2002. Supply chain management: a supplier perspective. *Production Planning and Control*, 13(4), pp. 370-380.

Rosetta, 2007. http://www.rosettanet.org, visited January 17, 2008.

Shirodkar S., Kempf, K., 2006. Supply Chain Collaboration Through shared Capacity Models, *Interfaces*, 36(5), pp. 420-432.

Småros, J., 2005. *Information sharing and collaborative forecasting in retail supply chains*, PhD thesis, Helsinki University of Technology, Laboratory of Industrial Management.

Tang, C.S., 2006. Perspectives in supply chain risk management, *International Journal of Production Economics*, 103, pp. 451-488.

Teixidor, A.B., 2006. *Proactive management of uncertainty to improve scheduling robustness in process industry*. Thèse de doctorat, Universitat Politèctina de Barcelone.

Thierry, C., Lauras, M., Lamothe, J., Mahmoudi, J., Charrel, J-C., 2006. Viewpoint-centred methodology for designing cooperation policies within a supply chain. *Proceedings from the International Conference on Information Systems, Logistics and Supply Chain (ILS 2006),* Botta-Genoulaz, V., Riane, F., (Eds.), Lyon, France.

Vics, 2007. http://www.vics.org/committees/cpfr/, visited January 17, 2008.

Ziegenbein, A., Nienhaus, J., 2004. Coping with supply chain risks on strategic, tactical and operational level. *Proceedings of the Global Project and Manufacturing Management Symposium*. Richard Harvey, Joana Geraldi, and Gerald Adlbrecht (Eds.), Siegen, pp. 165-180.